\documentclass[12pt]{amsart}

\usepackage{amsmath}
\usepackage{amssymb}

\def\Gal{\mathrm{Gal}}

\begin{document}

\title{Unramified extensions and geometric $\mathbb{Z}_p$-extensions 
of global function fields}
\author{Tsuyoshi Itoh}
\address{College of Science and Engineering,
Ritsumeikan University, 
1-1-1 Noji Higashi, Kusatsu, 
Shiga, 525-8577, JAPAN}
\email{tsitoh@se.ritsumei.ac.jp}
\maketitle

\begin{abstract}
We study on finite unramified extensions of
global function fields (function fields of one valuable over a finite field).
We show two results.
One is an extension of Perret's result about the ideal class group problem.
Another is a construction of a geometric $\mathbb{Z}_p$-extension which 
has a certain property.
\end{abstract}

\section{Main theorems}
Throughout the present paper, we fix a prime number $p$
and a finite field $\mathbb{F}$ of characteristic $p$.

It is known that there is a finite abelian group $G$ which 
does not appear as the divisor class group of degree 0 of 
global function fields (Stichtenoth \cite{Sti}).
On the other hand, Perret \cite{PerC} showed the following:

\par \bigskip

\noindent \textbf{Theorem A}(\cite{PerC}). {\itshape
For any given finite abelian group $G$, there is a finite separable geometric extension
$k/ \mathbb{F} (T)$ such that $\mathrm{Cl} (\mathcal{O}) \cong G$,
where $\mathcal{O}$ denotes the integral closure of $\mathbb{F} [T]$
in $k$ and $\mathrm{Cl}(\mathcal{O})$ the ideal class group of $\mathcal{O}$.
}

\par \bigskip

This theorem is shown by using the following:

\par \bigskip

\noindent \textbf{Theorem B}(\cite{PerC}). {\itshape
For any given finite abelian group $G$, there is a global function 
field $k$ over $\mathbb{F}$ and a finite set $S$ of places of $k$ such that $\mathrm{Cl}_S (k) \cong G$,
where $\mathrm{Cl}_S (k)$ denotes the $S$-class group of $k$.
}

\par \bigskip

Let $H_S (k)$ be the $S$-Hilbert class field of $k$, that is, 
the maximal unramified abelian extension field of $k$ in which 
all places of $S$ split completely (see \cite{RosH}).
We note that $\mathrm{Cl}_S (k) \cong \Gal (H_S (k)/k)$ 
by class field theory.
Hence Theorem B also implies the existence of $k$ and $S$ 
which satisfy $\Gal (H_S (k)/k) \cong G$.

In the present paper, we extend the above result for non-abelian cases.
We will show the following:

\par \bigskip

\noindent \textbf{Theorem 1. } {\itshape
For any given finite group $G$, 
there is a global function field $k$ over $\mathbb{F}$
and a finite set $S$ of places of $k$
such that $\Gal (\tilde{H}_S (k) /k) \cong G$,
where $\tilde{H}_S (k)$ denotes the maximal unramified
extension field over $k$ in which all places of $S$ split completely.
}

\par \bigskip

See also Ozaki \cite{OzaU} for the number field case.

We will prove Theorem 1 in section 2.
Our proof dues to Perret's idea (see \cite{PerC}).
That is, we will construct an unramified $G$-extension, and 
take a sufficiently large set $S$ of places 
such that $\Gal (\tilde{H}_S (k) /k) \cong G$.
(We use the term ``$G$-extension'' as a Galois extension
whose Galois group is isomorphic to $G$.)
To construct an unramified $G$-extension, we shall show an analogue of 
Fr{\"o}hlich's classical result \cite{Fro} for number fields.

In section 3, we shall apply Perret's idea to the Iwasawa theory.
Let $k$ be a global function field over $\mathbb{F}$,
$S$ a finite set of places in $k$, 
and $k_{\infty}/k$ a geometric $\mathbb{Z}_p$-extension.
(Recall that $p$ is the characteristic of $\mathbb{F}$.)
We assume that 
\begin{itemize}
\item[(A)] only finitely places of $k$ ramify in $k_{\infty}/k$, and 
\item[(B)] all places of $S$ split completely in $k_{\infty}/k$.
\end{itemize}
Under these assumptions,
we can treat the Iwasawa theory for the $S$-class group (see \cite{RosH}).
For a non-negative integer $n$, let $k_n$ be the $n$th layer of $k_{\infty}/k$.
That is, 
$k_n$ is the unique subfield of $k_\infty$ which is a cyclic extension over $k$ of degree $p^n$.
Moreover, let $A_n$ be the Sylow $p$-subgroup of the $S$-class group of $k_n$.
(Here we use the same character $S$ as the set of places of $k_n$ lying above $S$.)
We put $X_S = \varprojlim{A_n}$, and
we call $X_S$ the Iwasawa module of $k_{\infty}/k$ for the 
$S$-class group.
We put $\Lambda=\mathbb{Z}_p [[T]]$.
It is known that $X$ is a finitely generated torsion 
$\Lambda$-module, and the ``Iwasawa type formula''
holds for $A_n$ (see \cite{RosH}).
That is, there are non-negative integers $\lambda, \mu$, and an integer $\nu$ 
such that $|A_n| = p^{\lambda n + \mu p^n +\nu}$ for all sufficiently large $n$.

There is a natural problem: 
characterise the $\Lambda$-modules which appear as $X_S$.
(For the number field case, the same problem is dealt in, e.g.,  
\cite{OzaC}, \cite{FOO}.)
Concerning this problem, we shall give the following result including ``non-abelian'' cases.

\par \bigskip

\noindent \textbf{Theorem 2. } {\itshape
For any given finite $p$-group $G$, 
there exists a global function field $k$ over $\mathbb{F}$,
a finite set $S$ of places of $k$,
and a geometric $\mathbb{Z}_p$-extension $k_{\infty}/k$
such that $\Gal (\tilde{L}_S (k_n) /k_n) \cong G$ 
(as groups) for all $n \geq 0$,
where $\tilde{L}_S (k_n)$ is the maximal unramified 
pro-$p$-extension field over $k_n$ in which all places lying above $S$ split completely.
}

\par \bigskip

For the number field case, Ozaki \cite{OzaC} showed that 
every ``finite $\Lambda$-module'' appears as the Iwasawa module
of a $\mathbb{Z}_p$-extension.
In Theorem 2, if we take a finite abelian $p$-group as $G$, this is a weak analogue 
of Ozaki's result. 
That is, every finite $\Lambda$-module on which
$\Lambda$ acts trivially appears as $X_S$.

\section{Proof of Theorem 1}

\subsection{Function field analogue of Fr{\"o}hlich's Theorem}
At first, we shall show that for any finite group $G$,
there is an unramified geometric extension $K/k$ of global function fields
such that $\Gal(K/k) \cong G$.
For the number field case, Fr{\"o}hlich already showed the following result.
\par \bigskip

\noindent \textbf{Fr{\"o}hlich's Theorem}(\cite{Fro}). {\itshape
For every positive integer $n$, there is an unramified extension $K/k$ of algebraic number fields
such that $\Gal(K/k) \cong \mathfrak{S}_n$,
where $\mathfrak{S}_n$ denotes the symmetric group of degree $n$.
}

\par \bigskip

We will show the following:

\par \bigskip

\noindent \textbf{Theorem 3.} {\itshape
For every integer $n \geq 5$,
there is a global function field $k$ over $\mathbb{F}$
and an unramified geometric extension $K/k$
such that $\Gal(K/k) \cong \mathfrak{S}_n$.
}

\par \bigskip

To prove this, we follow Fr{\"o}hlich's original argument
(see also Malinin \cite{Mal}).
That is, we construct a certain $\mathfrak{S}_n$-extension over 
the rational function field $\mathbb{F}(T)$ and then 
we lift up this extension.

\par \bigskip

\noindent \textbf{Lemma 4. } {\itshape
Assume that $n \geq 5$.
There is a Galois extension $k'$ over $\mathbb{F}(T)$
which satisfies all of the following properties.
\begin{itemize}
\item $k'/\mathbb{F}(T)$ is an geometric extension.
\item $\Gal (k'/\mathbb{F}(T)) \cong \mathfrak{S}_n$.
\item $1/T$ is unramified in $k'/\mathbb{F}(T)$.
\end{itemize}
}

\par \bigskip

\noindent \textbf{Proof. } 
We first see that there is a $\mathfrak{S}_n$-extension over $\mathbb{F}(T)$.
This follows from the fact that $\mathbb{F}(T)$ is a Hilbertian field
(see, e.g, \cite[Corollary 16.2.7]{F-J}).

We put $A=\mathbb{F} [T]$.
Fix a monic separable polynomial $F(X) \in A[X]$ of degree $n$
such that the splitting field of $F(X)$
over $\mathbb{F}(T)$ is an $\mathfrak{S}_n$-extension.
We know that there is an element $N_F \in A$ which satisfies the following property: 
if a monic polynomial $G(X) \in A[X]$ of degree $n$ satisfies
$G(X) \equiv F(X) \pmod{N_F}$, then the splitting field of 
$G(X)$ over $\mathbb{F}(T)$ is also an $\mathfrak{S}_n$-extension.
Moreover, we can take $N_F$ which is prime to $T$.
We also fix such $N_F$.

To construct a geometric $\mathfrak{S}_n$-extension, we take $G(X)$ as 
follows:
\[ \begin{array}{ll}
G(X) \equiv F(X) & \pmod{N_F}, \\
G(X) \equiv \mbox{(distinct polynomials of degree $1$)} & \pmod{r},
\end{array} \]
where $r$ is a monic irreducible polynomial of $A=\mathbb{F} [T]$
such that the degree of $r$ is odd and $r$ is prime to $T N_F$.
By the first congruence, we see that the splitting field $k'$ of $G(X)$ 
is a $\mathfrak{S}_n$-extension.
We shall show that the coefficient field of $k'$ is $\mathbb{F}$.
Let $\overline{\mathbb{F}}$ be the algebraic closure of $\mathbb{F}$.
We note that $M:=k' \cap \overline{\mathbb{F}}(T)$ is a finite Galois 
extension over $\mathbb{F} (T)$.
Since $\Gal (k'/\mathbb{F} (T)) \cong \mathfrak{S}_n$ and $n \geq 5$, 
$M$ must be $\mathbb{F} (T)$ or the unique quadratic subfield in $k'/\mathbb{F} (T)$.
If $M \neq \mathbb{F} (T)$, then all places of odd degree do not split in $M$.
However, we see that the place generated by $r$
splits completely in $k'$ by the second congruence.
It is a contradiction. 

To satisfy the third condition, it is sufficient to show
that one can take $k'$ such that $T$ is unramified in $k'/\mathbb{F} (T)$.
(Because we replace an intermediate $T$ to $U=1/T$, then 
$1/U$ is unramified in $k'/\mathbb{F} (U)$ and the other conditions also satisfied.)
Then we take $G(X)$ as 
follows:
\[ \begin{array}{ll}
G(X) \equiv F(X) & \pmod{N_F}, \\
G(X) \equiv \mbox{(distinct polynomials of degree $1$)} & \pmod{r}, \\
G(X) \equiv \mbox{(an irreducible polynomial)} & \pmod{T}.
\end{array} \]
By the third congruence, we see that $T$ in unramified in $k'$. 
\hfill $\Box$

\par \bigskip

We shall prove Theorem 3.
We fix a geometric $\mathfrak{S}_n$-extension $k'/\mathbb{F} (T)$ satisfying the conditions given in Lemma 4.
Let $F(X) \in A[X]$ be the minimal polynomial of an generator of $k'$ over $\mathbb{F} (T)$.
$F(X)$ has degree $n!$ as a polynomial of $X$.

We define the following notations.
\begin{itemize}
\item $\{ \mathfrak{p}_1, \ldots, \mathfrak{p}_t \}$ : the set of places of $\mathbb{F} (T)$
which ramify in $k'$.
\item $\mathfrak{p}_{t+1}$ : a place of $\mathbb{F} (T)$ which is inert 
in the unique quadratic subextension of $k'/\mathbb{F} (T)$
(distinct from $\mathfrak{p}_1, \ldots, \mathfrak{p}_t$).
\item $\mathfrak{p}_{t+2}$ :  a place of $\mathbb{F} (T)$ which splits 
in the unique quadratic subextension of $k'/\mathbb{F} (T)$ and has \textbf{odd} degree
(distinct from $\mathfrak{p}_1, \ldots, \mathfrak{p}_t$).
\item $p_1, \ldots, p_{t+2}$ : irreducible monic polynomials of $A=\mathbb{F} [T]$ which generate
$\mathfrak{p}_1, \ldots, \mathfrak{p}_{t+2}$, respectively.
\end{itemize}
We shall give some remarks.
Since $1/T$ does not ramify in $k'$, we can take generators of above places 
as an element of $A$.
It is not trivial that one can really take $\mathfrak{p}_{t+1}, \mathfrak{p}_{t+2}$.
However, by using Theorem 9.13B of \cite{RosN} (which is a precise
version of the Chebotarev density theorem for global function fields),
we can take such places.

We put $m=n!$.
By using Lemma 4, we can also construct an $\mathfrak{S}_m$-extension 
over $\mathbb{F} (T)$.
Let $H(X)$ be a polynomial in $A[X]$ of degree $m$ 
which gives an $\mathfrak{S}_m$-extension.
Then there is an element $N_H$ of $A$ having the following property: 
if a monic polynomial $G(X) \in A[X]$ of degree $m$ satisfies
$G(X) \equiv H(X) \pmod{N_H}$, then the splitting field of 
$G(X)$ over $\mathbb{F}(T)$ is also an $\mathfrak{S}_m$-extension.
We can also take $N_H$ such that it is prime to $p_1, \ldots, p_{t+2}$.

We take a polynomial $G(X)$ of $A[X]$ (having degree $m$) which satisfy the
following conditions (1)--(4).
\begin{equation}
G(X) \equiv H(X) \pmod{N_H}.
\end{equation}
If $G(X)$ satisfies (1), then $G(X)$ gives a $\mathfrak{S}_m$-extension.
Let $L$ be the splitting field of $G(X)$ over $\mathbb{F} (T)$.
\begin{equation}
G(X) \equiv \mbox{(distinct polynomials of degree $1$)} \pmod{p_{t+1}}.
\end{equation}
If $G(X)$ satisfies (2), then we see that $\mathfrak{p}_{t+1}$ splits in the unique 
quadratic subfield of $L/\mathbb{F} (T)$.
On the other hand, $\mathfrak{p}_{t+1}$ is inert 
in the unique quadratic subextension of $k' /\mathbb{F} (T)$.
Since $\Gal (k' /\mathbb{F} (T)) \cong \mathfrak{S}_n$ and
$\Gal (L/\mathbb{F} (T)) \cong \mathfrak{S}_m$, we can see that 
$k' \cap L = \mathbb{F} (T)$, and then 
$\Gal (Lk'/L) \cong \mathfrak{S}_n$.
\begin{equation}
G(X) \equiv \mbox{(distinct polynomials of degree $1$)} \pmod{p_{t+2}}.
\end{equation}
If $G(X)$ satisfies (3), then the odd degree place $\mathfrak{p}_{t+2}$
splits completely in $Lk'/\mathbb{F} (T)$.
This implies that $Lk'/\mathbb{F} (T)$ is a geometric extension.
Finally, it is known that there is a positive integer $s_i$ for each $i=1, \ldots, t$
depending only on $F(X)$ such that if $G(X) \equiv F(X) \pmod{p_i^{s_i}}$ then
$L \mathbb{F} (T)_{\mathfrak{p}_i} = k' \mathbb{F} (T)_{\mathfrak{p}_i}$,
where $\mathbb{F} (T)_{\mathfrak{p}_i}$ is the 
completion of $\mathbb{F} (T)$ at $\mathfrak{p}_i$.
Hence if we take $G(X)$ satisfying 
\begin{equation}
G(X) \equiv F(X) \pmod{p_i^{s_i}} \;\; \mbox{for $i=1, \ldots, t$},
\end{equation}
then we can see that $Lk'/L$ is unramified at all places.

We can take $G(X)$ satisfying (1)--(4).
By the above arguments, the extension $Lk'/k'$ satisfies the assertion of 
Theorem 3.
\hfill $\Box$

\par \bigskip

\noindent \textbf{Remark. }
When $G$ is abelian, an unramified geometric $G$-extension 
was constructed by Angles \cite{Ang}.
Moret-Bailly \cite{Mor} also gives a result which is close to ours.

\par \bigskip

\subsection{Proof of Theorem 1}
By Theorem 3, we can construct an unramified extension with any given 
finite group as its Galois group.
Let $K/k$ be a geometric Galois unramified extension such that 
$\Gal (K/k) \cong G$.

\par \bigskip

\noindent \textbf{Proposition 5. } {\itshape
There is a finite set of places $S$ of $k$ such that 
(i) all places in $S$ split completely in $K$, and 
(ii) $\tilde{H}_S (k)/k$ is a finite extension.
}

\par \bigskip

\noindent \textbf{Proof. }
The crucial point of this proposition is choosing a set $S$ to satisfy (ii).
For a positive integer $N$, we put 
\[ B_N = \{ \mathfrak{p} \; \vert \; \mbox{$\mathfrak{p}$ is a place of $k$}, 
\deg(\mathfrak{p})=N, \mbox{$\mathfrak{p}$ splits completely in $K/k$.} \}. \]
Let $g$ be the genus of $k$, and $q$ the number of elements in $\mathbb{F}$.
If $N$ is sufficiently large, then we can see
\[ |B_N| > \frac{q^{N/2}-1}{N} \mathrm{Max} (g-1, 0) \]
by using Theorem 9.13B of \cite{RosN}.
We fix an integer $N$ which satisfies the above inequality.
According to Ihara's theorem \cite[Theorem 1(FF)]{Iha}, 
if $S \supset B_N$, then $\tilde{H}_S (k)/k$ is a finite extension.
Hence we can take $S$ to satisfy the conditions (i) and (ii).
\hfill $\Box$

\par \bigskip

The rest part of the proof of Theorem 1 is quite similar to Perret's
argument given in \cite{PerC}.
We choose a set $S$ of places which satisfies the conditions in 
Proposition 5.
For a nontrivial element $\sigma$ of $\Gal (\tilde{H}_S (k)/K)$,
we can take a place $\mathfrak{P}$ of $\tilde{H}_S (k)$ corresponding to $\sigma$
by the Chebotarev density theorem.
We can take $\mathfrak{P}$ which is unramified in $\tilde{H}_S (K)/K$.
Let $\mathfrak{p}$ be the place in $k$ which is lying below $\mathfrak{P}$.
Since the decomposition field 
of $\mathfrak{P}$ in $\tilde{H}_S (k)/k$ contains $K$ and $K/k$ is a Galois extension,
we see that $\mathfrak{p}$ splits completely in $K/k$.
Then we see $\tilde{H}_S (k) \supsetneq \tilde{H}_{S \cup \{\mathfrak{p}\}} (k) \supset K$.
Replacing $S\cup \{\mathfrak{p}\}$ to $S$ and repeating the above 
operation, we can obtain Theorem 1.
\hfill $\Box$

\par \bigskip

\noindent \textbf{Remark. }
Our construction also gives the fact that $\tilde{H}_S (k)/k$ is a geometric extension.

\section{Proof of Theorem 2}
Firstly, we shall show the following:

\par \bigskip

\noindent \textbf{Theorem 6. } {\itshape
Let $k$ be a finite Galois extension of $\mathbb{F} (T)$.
Then, there exists a finite set $S$ of places of $k$ and 
a geometric $\mathbb{Z}_p$-extension $k_{\infty}/k$ (which satisfies the
assumptions (A) and (B) in section 1) such that 
the Iwasawa module for the $S$-class group is trivial
(i.e., $\lambda=\mu=\nu=0$).
}

\par \bigskip

Precisely, we will show a slightly stronger result.
That is, we can take $k_\infty /k$ being the ``lift up'' of a geometric 
$\mathbb{Z}_p$-extension of $\mathbb{F} (T)$.
This fact is used to prove Theorem 2.

\par \bigskip

\noindent \textbf{Proof of Theorem 6.} 
We take a place $\mathfrak{p}_0$ of $\mathbb{F} (T)$ which splits completely in $k$.
We also take a place $\mathfrak{r}$ of $\mathbb{F} (T)$ which is distinct from $\mathfrak{p}_0$
and unramified in $k$.
We claim that there is a geometric $\mathbb{Z}_p$-extension $F_\infty/\mathbb{F} (T)$
unramified outside $\mathfrak{r}$ which satisfies that
\begin{itemize}
\item $\mathfrak{r}$ is totally ramified, and
\item $\mathfrak{p}_0$ splits completely.
\end{itemize}
We shall show this claim.
Let $M$ be the maximal pro-$p$-extension over $\mathbb{F} (T)$ which is unramified outside $\mathfrak{r}$.
Then we know that $\Gal (M/\mathbb{F} (T)) \cong \mathbb{Z}_p^{\infty}$
(see, e.g., \cite{G-K}).
Hence there are infinitely many geometric $\mathbb{Z}_p$-extensions 
which satisfy the above conditions.

Let $F_1$ be the initial layer of $F_\infty/\mathbb{F} (T)$, and 
we put $k_1=k F_1$.
Then $k_1/\mathbb{F} (T)$ is a Galois extension, and $\mathfrak{p}_0$
splits completely in $k_1$.
We set $S_0=\{ \mathfrak{p}_0 \}$, and we use the same character to denote 
the set of places lying above $\mathfrak{p}_0$.
We take an nontrivial element $\sigma_1$ of $\Gal (H_{S_0} (k_1)/k_1)$.

By using the above argument, we can take a geometric 
$\mathbb{Z}_p$-extension $F'_\infty/\mathbb{F} (T)$
unramified outside $\mathfrak{r}$ which satisfies
\begin{itemize}
\item $F'_\infty \cap F_\infty =\mathbb{F} (T)$,
\item $\mathfrak{r}$ is totally ramified in $F'_\infty F_\infty$, and
\item $\mathfrak{p}_0$ splits completely in $F'_\infty$.
\end{itemize}

Let $F'_1$ be the initial layer of $F'_\infty/\mathbb{F} (T)$.
Then we see that $F'_1 \cap k_1 = \mathbb{F} (T)$ and 
$k_1 F'_1 \cap H_{S_0} (k_1)=k_1$.
Let $\tau$ be a generator of the cyclic group $\Gal (F'_1/\mathbb{F} (T))$,
and $\tau_1$ an element of $\Gal (F'_1 H_{S_0} (k_1)/k_1)$
which is the image of $(\tau, \sigma_1)$ of the natural isomorphism
\[ \Gal (F'_1/\mathbb{F} (T)) \times \Gal (H_{S_0} (k_1)/k_1)
\longrightarrow \Gal (F'_1 H_{S_0} (k_1)/k_1). \]
We can regard $\tau$ as an element of $\Gal (F'_1 H_{S_0} (k_1)/\mathbb{F} (T))$.
By the Chebotarev density theorem, there is a place $\mathfrak{P}_1$ of 
$F'_1 H_{S_0} (k_1)$ which corresponds to $\tau_1$.
Let $\mathfrak{p}_1$ be the place of $\mathbb{F} (T)$ lying below $\mathfrak{P}_1$.
Then we see that $\mathfrak{p}_1$ splits completely in $k_1$ and is inert in $F'_1$.
We put $S_1 = S \cup \{ \mathfrak{p}_1 \}$.

We do not know whether $\mathfrak{p}_1$ splits completely in $F_\infty$ or not.
It is a problem because we need the assumption (B) in section 1.
To evade this problem, 
we replace $F_\infty$ to another geometric $\mathbb{Z}_p$-extension.
We remark that $F_\infty F'_\infty/\mathbb{F} (T)$ is a $\mathbb{Z}_p^2$-extension
unramified outside $\mathfrak{r}$.
Since $\mathfrak{p}_1$ does not split in $F'_1$, it also does not split in $F'_\infty$.
Hence the decomposition field of $F_\infty F'_\infty/\mathbb{F} (T)$ for $\mathfrak{p}_1$ is 
a $\mathbb{Z}_p$-extension over $\mathbb{F} (T)$.
We denote it $F''_\infty$.
We also note that $F''_\infty/\mathbb{F}(T)$ is the unique $\mathbb{Z}_p$-extension 
contained in $F_\infty F'_\infty$ such that $\mathfrak{p}_1$ splits completely.
Then the initial layer of $F''_\infty/\mathbb{F}(T)$ must coincide with $F_1$.
We replace $F_\infty$ to $F''_\infty$.

We note that $H_{S_0} (k_1) \supsetneq H_{S_1} (k_1)$
by the definition of $\mathfrak{p}_1$.
Similarly, we can choose a place $\mathfrak{p}_2$, put $S_2 = S_1 \cup \{ \mathfrak{p}_2 \}$, and replace 
a $\mathbb{Z}_p$-extension such that all places in $S_2$ splits completely.
Repeating this operation, we see that $H_{S_t} (k_1) =k_1$ for some finite set $S_t$.
We note that $F_{\infty} k/k$ satisfies the assumptions (A) and (B).

Finally, we shall give an Iwasawa-theoretic argument.
In $F_\infty k/k$, all ramified places (these are lying above $\mathfrak{r}$) are totally ramified.
From this, we also see $H_{S_t} (k) =k$.
Let $k_n$ be the $n$th layer of $F_\infty k/k$, and
$A_n$ the Sylow $p$-subgroup of $\mathrm{Cl}_{S_t} (k_n)$.
By the above results, we see that both of $A_0$ and $A_1$ are trivial.
In this situation, we can use the method given by Fukuda \cite{FukR}.
Hence we can obtain the fact that $A_n$ is trivial for all $n$.
This implies the assertion of Theorem 6.
\hfill $\Box$

\par \bigskip

We shall show Theorem 2.
We fix a finite $p$-group $G$.
From the proof of Theorem 1, we can take a Galois extension $K/\mathbb{F} (T)$
and a subfield $k$ of $K$ such that $K/k$ is unramified and 
$\Gal (K/k) \cong G$.
From the proof of Theorem 6, we can take a geometric $\mathbb{Z}_p$-extension
$F_\infty/ \mathbb{F} (T)$ such that $F_\infty \cap K = \mathbb{F} (T)$, 
and a set $S$ of places (of $\mathbb{F} (T)$) such that 
the order of the $S$-class group of every layer of $F_\infty K/K$ is prime to $p$.
Since the $p$-group $G$ is solvable, the $\mathbb{Z}_p$-extension
$F_\infty k/k$ satisfies the assertion of Theorem 2.
\hfill $\Box$

{}

\end{document}